\documentclass[11pt]{article}
\usepackage{amsfonts}

\newtheorem{Theorem}{Theorem}

\newtheorem{Lemma}[Theorem]{Lemma}
\newtheorem{Proposition}[Theorem]{Proposition}

\newcommand{\R}{\mathbb{R}}

\newcommand{\Z}{\mathbb{Z}}
\newcommand{\E}{\mathbb{E}}
\newcommand{\Pp}{\mathbb{P}}

\newcommand{\Eta}{\eta}
\newcommand{\erf}{\operatorname{erf}}
\newcommand{\erfc}{\operatorname{erfc}}
\newcommand{\K}{\mathbf{K}}
\newcommand{\J}{\mathbf{J}}

\newcommand{\Pf}{\operatorname{Pf}}
\newcommand{\qed}{\vrule height7.5pt width4.17pt depth0pt}

\newcommand{\phfac}{\phi}
\newcommand{\Phfac}{\Phi}

\newcommand{\ind}{\mathbf{I}}

\usepackage{amsmath} %for substack
\usepackage{hyperref} %for clickable links

\begin{document}

\title{Examples of interacting particle systems on $\Z$ as Pfaffian point processes: coalescing
 branching random walks and annihilating random walks with immigration.}
 
 \author{
 Barnaby Garrod, Roger Tribe, Oleg Zaboronski \\ \\
 Mathematics Institute, University of Warwick}

\date{}

\maketitle

\begin{abstract}
Two classes of interacting particle systems on $\Z$ are shown to be Pfaffian point processes, at any fixed time and 
for all deterministic initial conditions. The first comprises coalescing and branching random walks, 
the second annihilating random walks with pairwise immigration. Various limiting Pfaffian point processes on 
$\R$ are found by diffusive rescaling, including the point set process for the Brownian web and Brownian net. 
%For certain maximal initial conditions the processes may be related by thickening and thinning.
\end{abstract}

\noindent

\section{Introduction and statement of key results}  \label{s1}
Interacting particle systems of reaction-diffusion type have been an object of mathematical investigation for a long time,
see e. g. \cite{griffeath}. It turns out that many such systems can be studied using methods of integrable probability.
%added the introductory sentence
This paper continues the study in \cite{GPTZ:15}, where it was shown that systems of instantaneously annihilating or coalescing random walks on
$\Z$ are Pfaffian point processes, at any fixed time and for all deterministic initial conditions. 
%The proofs involves finding a suitable Markov duality function, 
%and exploiting this to calculate the particle intensities as Pfaffians.  Additional natural mechanisms, such as allowing particles to die at 
%constant rates, or introducing a Poisson immigration of particles, seem to destroy the Pfaffain point process property. 
The purpose in this paper is to describe two additional mechanisms that preserve this Pfaffian property.  
The Pfaffian property should be useful to investigate statistics, such as asymptotics of correlation functions 
(as in Theorem 1 of \cite{TZ:11}) or for studying gap probabilities (as in \cite{TW}, \cite{RS} or \cite{TZ3} for 
examples from random matrix ensembles).  
 
The algebraic structure of the generators of various one dimensional interacting particle systems 
(without necessarily preservation of particle numbers) has been investigated, and is in many examples intimately 
linked to Hecke algebras, see the reviews \cite{chemsurvey2} and \cite{chemsurvey1}.  
A future aim is to better understand how to deduce from these algebraic properties 
the concrete statistical properties of our models, such as the correlation functions, and in particular what algebraic 
properties lie behind the emergence of the Pfaffian property.

The connection between coalescing and annihilating systems is well known, 
both at the analytic level and
also via direct couplings of the two processes, see \cite{swart} for review. 
%added reference to Swart
Such pathwise couplings 
are fundamental to the study of the 
Brownian Web and Net and the dual Brownian Web and Net (and their discrete analogues - see \cite{SS:08} and \cite{SSS:15}).
Moreover all the models in this paper can be found 'within' the Net, and so there
will be many pathwise couplings available. Our results can be seen as analytic consequences of
a strong coupling between our two models, namely coalescing plus branching random walks 
and annihilating random walks plus pairwise immigration, with finite systems of annihilating random walks. 
We do not exploit this coupling, but we concentrate on investigating 'fixed time' laws, aimed at revealing an 
algebraic structure behind the processes.\footnote{A colouring construction of annihilating particles starting
from a coalescing system can be used to relate the Pfaffian point processes describing the fixed time laws of coalescing
and annihilating random walks or Brownian motions, see \cite{TZ:11}. At the moment we do not see a similar link
between branching-coalescing random walks and annihilating random walks with pairwise immigration, but a further
investigation might uncover such a link.}
%added a footnote

\subsection{(BCRW) Branching coalescing random walks.} \label{s1.1}
%\textbf{(BCRW) Branching coalescing random walks.}  
In addition to instantaneously coalescing random walks, we allow nearest neighbour binary branching: any given particle may instantaneously produce an independent copy at a nearest neighbour. The dynamics of this continuous-time model on $\Z$ are informally described as follows. Between interactions particles perform independent nearest neighbour random walks with jumps
\[
 \mbox{$x \to x-1$ at rate $q$, \hspace{.2in} and  \hspace{.2in}  $x \to x+1$ at rate $p$.}
\]
If a particle jumps onto an occupied site then the two particles instantaneously coalesce. Independently a particle branches
\[
 \mbox{$x \to \{x,x-1\}$ at rate $\ell$, \hspace{.2in} and  \hspace{.2in} $x \to \{x,x+1\}$ at rate $r$.}
\]
Branching events respect coalescence: if a particle branches onto an occupied site then the existing and new particles instantaneously coalesce.

The generator is given in (\ref{eq:processgeneratorbcrw}), and this characterises the law of a process with values in 
$\{0,1\}^{\Z}$ which we denote as the BCRW model.  We write $(\eta_t:t \geq 0)$ for canonical variables and $\Pp_{\eta}$ for the law corresponding to an initial 
condition $\eta \in \{0,1\}^{\Z}$. As in  \cite{GPTZ:15}, we fix $t \geq 0$ and consider $\eta_t$ as a point process on $\Z$, and our aim is to 
establish that this is a Pfaffian point process for a suitable Pfaffian kernel $\K(x,y)$. 
 A simple point process on $\Z$ is called Pfaffian if there is a $2\times 2$ complex matrix-valued kernel $\K$
on $\Z \times \Z$, $(x,y) \mapsto (K_{\alpha \beta}(x,y))_{1\leq \alpha \beta\leq 2}$, such that the $n$-th
intensity function of the process $\rho^{(n)}$ is given by 
$\rho^{(n)}(x_1, x_2, \ldots, x_n)=\Pf[\K(x_i,x_j)]_{1\leq i,j \leq n}, n=1,2,\ldots$
(see \cite{zeitouni} for more details). 

Recall the notation from \cite{GPTZ:15}: for $\eta \in \{0,1\}^{\Z}$ and $y \leq z$ we write
\[
\eta[y,z) = \sum_{y \leq x<z} \eta(x) 
\]
so that $\eta[y,y) = 0$. For $f:\Z \to \R$ we write difference operators as
\[
D^+f(x) = f(x+1)-f(x), \quad D^- f(x) = f(x-1)-f(x). 
\]
If the difference operator $D^{\pm}$ is applied to the $i$-th argument of a 
function of several variables we use the notation $D^{\pm}_{i}$.
\begin{Theorem} \label{theorem:BCRW}
Let $(\eta_t:t \geq 0)$ be the BCRW model with parameter values satisfying
\begin{equation}\label{condition:rates}
 p\ell=qr,\qquad p,q>0.
\end{equation}
For any initial condition $\eta\in\{0,1\}^{\Z}$, and at any fixed time $t \geq 0$, the variable $\eta_t$ is a Pfaffian point process with kernel $\K$ given, for $y<z$, by
 \begin{equation} \label{discretekernel1}
\K(y,z)
  = \frac{1}{\phfac}\left(\begin{array}{cc}
           K_t(y,z) & -D^+_2 K_t(y,z) \\
           - D^+_1 K_t(y,z) & D^+_1 D^+_2 K_t(y,z)
           \end{array}\right),
 \end{equation}
and $
\K_{12}(y,y)=1 - \frac{1}{\phi} K_t(y,y+1)$, with all other entries determined by anti-symmetry, where 
\[
K_t(y,z) = \phi^{z-y} \Pp_{\eta}[ \eta_t[y,z)=0], \quad \mbox{for $t \geq 0$ and $y < z$,}
\]
and $ \phi = \sqrt{1+ \frac{\ell}{q}} = \sqrt{1+\frac{r}{p}}$. 
The same result holds when $p=r=0$ and $q,l>0$ by taking
$ \phi = \sqrt{1+ \frac{\ell}{q}}$, or when  $q=l=0$ and $p,r>0$ by taking
$ \phi = \sqrt{1+\frac{r}{p}}$.
\end{Theorem}
\textbf{Remarks.}  \textbf{1.1.} 
For certain random initial conditions, including the natural case when the sites $(\eta_0(x):x \in \Z)$ are independent,
the process does remain a Pfaffian point process  (see the remarks after the proof of Lemma \ref{lemma:emptyintervalpfaffian}). The invariant measure, which is product Bernouilli, can be considered as a Pfaffian point process 
(see the remark at the end of section \ref{s2}.)

\noindent
\textbf{1.2.}  The restriction $pl=qr$ seems to be necessary for the Pfaffian point process property,
and we don't fully understand its origin. Luckily, 
in the large scale diffusive limits explained below this restriction plays no role. Indeed 
the limit continuum systems depend on three parameters. This is consistent with
the three parameters that can be used in the Brownian net, for example the
two drift parameters and a common diffusion parameter for the left
and right paths, which together
determine the branching rate, see \cite{SS:08}.
For the four parameter lattice models, there
exist Markov duality functions even when the the restriction $pl=qr$ is not true, and it would be of interest to 
examine whether some other algebraic structure holds for the correlation functions.

\noindent
\textbf{1.3.} 
In section {\ref{s4} we investigate continuum limits to our Pfaffian point kernels, under space time diffusive scaling, yielding Pfaffian point processes on $\R$.
For the branching parameter to have an effect in the continuum kernel, it needs to be scaled to grow suitably fast. 
This however is well understood in the construction of the Brownian net (see Sun and Swart \cite{SS:08}),
and we of course need the same scaling of the branching as in the discrete time branching random walk approximations to the net. The Brownian net is a continuum 
collection of space time paths found by scaling discrete time branching coalescing random walks started at all space-time lattice points. 
The point set process $(\xi^A_t:t \geq 0)$ within the net, is defined by 
looking at all points at time $t$ that are on paths that start at time zero in a set $A \subset \R$. It is known to be a Feller process 
taking values in the compact sets of $\overline{\R}$ with a suitable Hausdorff metric (see Theorem 1.11 in \cite{SS:08})). The approximating discrete time branching coalescing 
random walks used in the construction of the net are not Pfaffian point processes (the discrete time difference equations
analogous to the Markov duality below are not solved by Pfaffians). However the difference between the discrete time models and continuous time
models does not affect the continuum limits. We leave the verification of these technical details to a forthcoming paper, 
but we state here the resulting Pfaffian property for the Brownian net point set process, which answers the first 
open problem in section $8.3$ of the survey paper on the net \cite{SSS:15}. The Brownian net can have a
parameter $b \geq 0$, controlling the branching rate, where the embedded left-right paths have drifts $\pm b$. 
The standard Brownian net corresponds to $b=1$, and the Brownian web corresponds to the case $b=0$.
\begin{Proposition} \label{net}
 The transition density $p_t(A,dB)$ for the point set process $(\xi^A_t)$, within the Brownian net with parameter $b \geq 0$, 
is equal to the distribution of the closed support of a Pfaffian point process on $\overline{\R}$ 
with the kernel $\K^A_t$ of the form
 \begin{eqnarray} 
&& \K_t^{A}(y,z)
  = \left(\begin{array}{cc}
           K^{A}_t(y,z) & -D_2 K^{A}_t(y,z) \\
           - D_1 K^{A}_t(y,z) & D_1 D_2 K^{A}_t(y,z)
           \end{array}\right)
 \quad \mbox{for $y<z$,} \nonumber \\
 && (\K_t^{A})_{12}(y,y) = -  D_2 K^{A}_t(y,y) + b, \label{BNkernel}
 \end{eqnarray}
 where $D_i$ is a partial derivative in the $i$th variable, and  
where $(K^{A}_t(y,z):y,z\in \R^2:y \leq z)$ is the unique bounded solution to the PDE
\begin{equation}
\left\{ \begin{array}{rcl}
\partial_t K^{A}_t(y,z) &=&  \frac12 \Delta  K^{A}_t(y,z)  - b^2 K^{A}_t(y,z)\\
K^{A}_t(y,y) & = & 1 \\
K^{A}_0(y,z) & = & \ind((y,z) \cap A = \emptyset)e^{b(z-y)}.
\end{array} \right. \label{netpde}
\end{equation}
\end{Proposition}
The examples in section \ref{s4} will illustrate how continuum kernels arise of this form. 

\subsection{(ARWPI) Annihilating random walks with pairwise immigration.} \label{s1.2}
As for the previous model, particles jump left or right at rates $q$ and $p$. 
If a particle jumps onto an occupied site then the two particles instantaneously annihilate. 
In addition, independently there is
\[
 \mbox{immigration of a pair of particles on sites $\{x,x-1\}$ at rate $m$.}
\]
Immigration respects annihilation: if a particle immigrates onto an occupied site then the existing and new particles 
instantaneously annihilate. The generator is given in (\ref{eq:processgeneratorarwpi}) and characterises a Markov process that we denote as the ARWPI model. 
 
\begin{Theorem} \label{theorem:ARWPI}
 For any initial condition $\eta\in\{0,1\}^{\Z}$ for the ARPI model, and at any fixed time $t\geq 0$, 
 the variable $\eta_t$ 
 is a Pfaffian point process with kernel $\K$ given, for $y<z$, by
 \begin{equation} \label{discretekernel2}
  \K(y,z) = \frac{1}{2}\left(\begin{array}{cc}
	      K_t(y,z) & - D^+_2K_t(y,z) \\
	      - D^+_1K_t(y,z) & D^+_1D^+_2K_t(y,z) 
	    \end{array}\right),
 \end{equation}
 and $\K_{12}(y,y)=-\frac{1}{2}D^+_2K_t(y,y)$, with other entries determined by anti-symmetry, where 
\[
K_t(y,z) = \E_{\eta}[ (-1)^{\eta_t[y,z)}], \quad \mbox{for $t \geq 0$ and $y < z$.}
\]
\end{Theorem}
\textbf{Remarks.}  \textbf{1.4.} The Glauber spin chain on $\Z$ is an assignment of $\pm1$ spin values to each site which independently
flip according to rates determined by nearest neighbour spins \cite{glauber}. Sites favour aligned spin and at zero temperature a site surrounded
by spins of the same sign does not flip, and the domain wall between regions of constant spin form a system of 
annihilating random walks on the dual lattice. At positive temperature, a spin may spontaneously
flip regardless of its neighbours, and this corresponds to the creation of a pair of neighbouring domain walls. 
Since the Glauber model can be solved at all temperatures by mapping to a system of free 
fermions (Felderhof \cite{felderhof}), it is 
reasonable that the extra immigration of pairs does not destroy the Pfaffian property of solutions. 
A model with Poisson immigration of single particles is perhaps of more interest, but we do not see a simple
algebraic structure behind this model.

\textbf{1.5.}  We give in section \ref{s3} a spatially inhomogeneous version of Theorem \ref{theorem:ARWPI}, where the parameters $p_x,q_x,m_x$ may be site dependent,
in particular allowing immigration of particles at different rates at different places. 
Continuum limits can also be found for the ARWPI model under diffusive rescaling, where the immigration rate $m$ must be scaled suitably so that they have a non-trivial effect in the limit. It is a pleasant fact that 
in many cases the Pfaffian kernel can be found completely explicitly. For example in example (d) in section \ref{s4}, 
which we call the Brownian firework,  pairs of particles are immigrated only at the origin at an infinite 
rate. This has a steady state $X^{(c)}_{\infty}$ where the immigration and the annihilation balance each other: 
$X^{(c)}_{\infty}$ is a Pfaffian point process 
on $\R \setminus \{0\}$ with kernel $\K^{(c)}_{\infty}$ of the form
 \begin{eqnarray*} 
&& \K_{\infty}^{(c)}(y,z)
  = \frac12 \left(\begin{array}{cc}
           K^{(c)}_{\infty}(y,z) & -D_2 K^{(c)}_{\infty}(y,z) \\
           - D_1 K^{(c)}_{\infty}(y,z) & D_1 D_2 K^{(c)}_{\infty}(y,z)
           \end{array}\right)
 \quad \mbox{for $y<z$,} \nonumber \\
 && (\K_{\infty}^{(c)})_{12}(y,y) = - \frac12 D_2 K^{(c)}_{\infty}(y,y) 
 \end{eqnarray*}
where
\begin{equation} \label{firework}
K^{(c)}_{\infty}(y,z) = \left\{ \begin{array}{ll}
1 + \frac{2}{\pi} \left( \arctan \frac{y}{z} - \arctan \frac{z}{y} \right) & \mbox{when $0<y<z$,} \\
0 & \mbox{when $y<0<z$,} \\
1 + \frac{2}{\pi} \left( \arctan \frac{z}{y} - \arctan \frac{y}{z} \right) & \mbox{when $y<z<0$,} \\
\end{array} \right.
\end{equation}
The corresponding intensity is given by $ \rho^{(1)}_{\infty}(y) = \frac{1}{\pi |y|}$.
% suggesting that there may be an accumulation point at the origin. 
Moreover since $\K^{(c)}_{\infty}(y,z) = 0$ when $y<0<z$ it is simple to deduce that
\[
\rho^{(n+m)}_{\infty} (y_1,\ldots,y_n,z_1,\ldots,z_m) = 
\rho^{(n)}_{\infty} (y_1,\ldots,y_n)
\rho^{(m)}_{\infty} (z_1,\ldots,z_m)
\] 
when $y_1,\ldots,y_n<0<z_1,\ldots,z_m$
and hence that $X^{(c)}_{\infty} |_{(-\infty,0)}$ and $X^{(c)}_{\infty} |_{(0,\infty)}$
are independent point processes. The infinite strength firework of particles at the origin 
leads to the two half spaces being independent. 

\textbf{1.6.} One simple consequence of the Pfaffian structure is an
estimate showing exponential convergence to equilibrium (which is a product Bernouilli distribution). 
Indeed,
writing $\rho^{(N)}_t( y_1,\ldots,y_N)$ for the $N$ point intensity function of the process as time $t$, 
we claim that there exist $C_N$ for all $N \geq 1$ so that
 \begin{equation} \label{expconv}
\left| \rho^{(N)}_t( y_1,\ldots,y_N) - \rho^{(N)}_{\infty}( y_1,\ldots,y_N) \right| \leq C_N e^{-2mt} 
\end{equation}
uniformly over for all $y_1,\ldots,y_N$ and over all initial conditions. (Recall that $m$ is the immigration rate of pairs). 
This follows for a deterministic initial condition
$\eta$ once we show that
\begin{equation} \label{kernelconv}
\left| K_t(y,z) - K_{\infty}(y,z)\right| \leq 2 e^{-2mt} \quad \mbox{for all $y,z \in \Z$,} 
\end{equation}
since the entries in the kernel $\K$ are differences of the bounded function $K_t$, so that the Pfaffian formula for
$\rho^{(N)}_t$ is given by a finite linear combination of finite products of $K_t(y_i,y_j)$. 
For a general initial condition, one can first condition on the initial condition
and then use the fact that the above estimates are uniform in $\eta$. 

(To show (\ref{kernelconv}) we can solve for $K_t(y,z)$ explicitly.
Indeed, fixing a deterministic initial condition $\eta$, the kernel $K_t(y,z)$ has a 
representation in terms of a pair of independent continuous time random walkers $(Y_t,Z_t)$ with 
generator $qD^+ + p D^-$, started at $Y_0=y, Z_0=z$. Let $\tau = \inf\{t:Y_t=Z_t\}$. Then
the equation solved by $K_t(y,z)$ (see Lemma \ref{lemma:MarkovdualityA}) shows that
\[
K_t(y,z) = \E \left[ e^{-2m \tau} \ind(\tau \leq t)\right] + e^{-2m t} \, \E \left[ K_0(Y_t,Z_t) \ind(\tau >t) \right]
\]
where $K_0(y,z) = (-1)^{\eta[y,z)}$ is bounded by $1$. The limit $K_{\infty}(y,z) = \E[ e^{-2m \tau}] $ 
and the estimate (\ref{kernelconv}) follows easily from subtracting these two probabilistic representations.
Solving explicitly we have 
\[
K_{\infty}(y,z) = \theta^{z-y}, \quad \mbox{where $\theta \in (0,1)$ solves $\theta+\theta^{-1} -2 = \frac{2m}{p+q}$.}
\]
The Pfaffian kernel of the form (\ref{discretekernel2}) corresponding to $K_{\infty}(y,z)$  is 
\[
  \K(y,z) = \frac{\theta^{z-y}}{2}\left(\begin{array}{cc}
	      1 & (1-\theta) \\
	    (1-\theta^{-1}) & (1-\theta)(1-\theta^{-1})
	    \end{array}\right),
\]
and $\K_{12}(y,y) = (1-\theta)/2$.
A little manipulation shows that this is a kernel for a product Bernouilli($\hat{\theta}$) distribution,
where 
\[
\hat{\theta} = \frac{1-\theta}{2} = \frac12 \left( \sqrt{\frac{m^2}{(p+q)^2} + \frac{2m}{p+q}} - \frac{m}{p+q} \right)
\] 
Indeed by conjugating $E^T \K E$
with an elementary matrix $E$ for row and column operations (which leaves the corresponding point process unaltered) 
$\K(y,z)$ for $y<z$ can be changed successively to
\begin{equation} \label{star1}
\K(y,z)  \to \frac{\theta^{z-y}}{2}\left(\begin{array}{cc}
	      1 & 0 \\
	        \theta-\theta^{-1}& 0
	    \end{array}\right) \to 
	    \frac{\theta^{z-y}}{2}\left(\begin{array}{cc}
	      0 & 0 \\
	      \theta-\theta^{-1} & 0
	    \end{array}\right),
\end{equation}
while leaving $\K_{12}(y,y) = (1-\theta)/2$ unchanged. Then the Pfaffian for $\rho^{(N)}$ has only a single non-zero entry
on the top row, and expanding along this row one finds 
\[
\rho^{(N)}(y_1,\ldots,y_N) = \frac{1-\theta}{2} \rho^{(N-1)}(y_2,\ldots,y_N).
\]
We remark that no exponential convergence statement such as (\ref{expconv}) holds for the BCRW model, since long 
empty gaps in the initial condition are only filled at linear speed. However, for many initial conditions there is weak convergence to a Bernoulli invariant measure - see the remarks in section \ref{s2}. 
%The argument breaks down since, although the discrete equation for the kernel is similar, 
%the initial conditions are no longer bounded.

%
%
\section{Proof of Theorem \ref{theorem:BCRW}.}  \label{s2}
We start with a terse summary of the main steps: the proof in \cite{GPTZ:15} for coalescing systems 
without branching uses the empty interval duality function;
this function remains a duality function for the branching model, but the empty interval probabilities are no longer 
given by a Pfaffian; however the duality function can be adjusted by a suitable phase factor in a way that again yields
 Pfaffians. The use of empty interval probabilities to study branching systems is not new, 
 see for example Krebs et al. \cite{krebs} where the equations for a single empty interval are used to study various finite systems.
 
The generator for BCRW is given, for suitable $F:\{0,1\}^{\Z}\to\R$, by
\begin{eqnarray}   
 \mathcal{L}  F(\eta) & = & q\sum_{x\in\Z} (F(\eta_{x,x-1})-F(\eta) ) + p\sum_{x\in\Z} (F(\eta_{x-1,x}) - F(\eta)) \nonumber \\
&& \hspace{.2in} + \ell\sum_{x\in\Z} (F(\eta^{b}_{x,x-1})-F(\eta) ) + r\sum_{x\in\Z} (F(\eta^{b}_{x,x+1}) - F(\eta)),
\label{eq:processgeneratorbcrw}
\end{eqnarray}
where $\eta_{x,y}$ (resp.\@ $\eta^{b}_{x,y}$) is the new configuration resulting from a jump (respectively a branch) from $x$ onto $y$, that is 
\[
\left\{
\begin{array}{l}
\eta_{x,y}(z) = \eta^{b}_{x,y}(z) = \eta(z) \quad \mbox{for $z \not \in \{x,y\}$,} \\ 
\eta_{x,y}(x) = 0, \quad
\eta^b_{x,y}(x) = \eta(x),\\
\eta_{x,y}(y) = \eta^b_{x,y}(y) = \min\{1, \eta(x)+\eta(y)\}.
\end{array} \right.
\]

For $n\ge1$ and $y=(y_1,\dots,y_{2n})$ with $y_1\le\dots\le y_{2n}$ we define the function $ \Sigma_y(\eta)$ as the inidcator that the intervals $[y_1,y_2), \ldots, [y_{2n-1},y_{2n})$ are all empty; explicitly
\[
 \Sigma_y(\eta) = \prod_{i=1}^n \ind \left(\eta[y_{2i-1},y_{2i})=0\right).
\]
 We define a one-particle operator, acting on $f:\Z\to\R$, by
\[
 L^{p,q} f(x) =  qD^+f(x) + pD^-f(x).
\]

\begin{Lemma}\label{lemma:Markovduality}
 For $y_1<\dots<y_{2n}$ the action of the generator $\mathcal{L}$ on $\Sigma_y(\eta)$ is given by
 \[
 \mathcal{L} \Sigma_y(\eta)
 = \sum_{i=1}^n (L^{p+r,q}_{y_{2i-1}} + L^{p,q+l}_{y_{2i}}) \, \Sigma_y(\eta)
\]
where the subscript $y_i$ indicates the variable upon which the operator acts.
\end{Lemma}

\textbf{Proof of Lemma~\ref{lemma:Markovduality}.}
A direct check shows that
the terms of $\mathcal{L}$ coming from left and right jumping contribute
\[
 \begin{array}{l}
  q \sum_{x\in\Z} \left (\Sigma_y(\eta_{x,x-1})-\Sigma_y(\eta)  \right) + p\sum_{x\in\Z} \left(\Sigma_y(\eta_{x-1,x}) - \Sigma_y(\eta)\right) 
  = \sum_{i=1}^{2n} L^{p,q}_{y_i}\Sigma_y(\eta)
 \end{array}
\]
to $  \mathcal{L}\Sigma_y(\eta)$ (see \cite{GPTZ:15} for the details of this calculation). It remains to compute the terms arising from branching. 
Consider the term from left branching. The modified branching configuration $\eta^b_{x,x-1}$ differs from $\eta$ only at the 
 site $x-1$, so for each $x$ there can be a change in at most one of the indicators in $\Sigma_y$, so we may write
\begin{eqnarray*} 
&&  \hspace{-.2in} \ell \left( \Sigma_y(\eta^{b}_{x,x-1})-\Sigma_y(\eta) \right) \\ \nonumber 
&=& \ell \sum_{i=1}^n \left(\prod_{j=1,j\neq i}^n \ind \left(\eta[y_{2j-1},y_{2j})=0\right) \right) 
\left (
\ind \left(\eta^b_{x,x-1}[y_{2i-1},y_{2i})=0\right) - \ind \left(\eta[y_{2i-1},y_{2i})=0\right)
\right).   %\label{eq:generator_separate}
\end{eqnarray*}
Fix $y<z$ and consider the generator contribution for a single empty interval indicator, namely
\[
 \ell \sum_{x\in\Z} \left(  \ind (\eta^b_{x,x-1}[y,z)=0) -  \ind \left(\eta[y,z)=0\right)  \right).
\]
The terms indexed by $x\le y$ and $x\ge z+1$ are zero, as the modified configuration is unchanged in the interval $[y,z)$.
 The terms indexed by $y\le x\le z-1$ are also zero since there must be a particle at $x$ to branch to the 
left from, in which case both empty interval indicators are zero. The remaining summand, when $x=z$, is given by
\begin{eqnarray*}
\ell \, \ind (\eta^{b}_{z,z-1}[y,z)=0 ) - \ell \, \ind \left(\eta[y,z)=0\right)
&=&  \ell \left( (1 - \eta(z)) - 1 \right) \ind \left(\eta[y,z)=0\right) \\
& = & \ell D^+_z \ind \left(\eta[y,z)=0\right).
\end{eqnarray*}
A similar calculation reveals that the term of the generator arising from right branching satisfies
\[
 r \sum_{x\in\Z} \left (
\ind (\eta^{b}_{x,x+1}[y,z)=0)- \ind \left(\eta[y,z)=0\right)  
\right)
 = r D^-_y \ind \left(\eta[y,z)=0\right).
\] 
Collecting up contributions gives the claimed action. \qed

The expression for $ \mathcal{L}\Sigma_y(\eta)$ in Lemma \ref{lemma:Markovduality} has different operators acting on even co-ordinates $y_{2i}$ and odd co-ordinates $y_{2i-1}$. The proof of the Pfaffian property in 
Lemma \ref{lemma:emptyintervalpfaffian} below
is facilitated if each co-ordinate has the same operator acting on it, and the aim is to introduce a suitable phase factor precisely to have this effect. The phase factor is defined by 
\[
\Phfac(y)=\prod_{i=1}^n\phfac^{(y_{2i}-y_{2i-1})} \quad \mbox{for $y=(y_1,\ldots,y_{2n}) \in \R^{2n}$ and $n \geq 1$}
\] 
and the following lemma shows that the value $ \phi = \sqrt{1+ \frac{\ell}{q}} = \sqrt{1+\frac{r}{p}}$ is the correct choice.
\begin{Lemma} \label{lemma:balancing}
Suppose the rates $p,q,r,l$ satisfy (\ref{condition:rates}) and $\phi$ is chosen as in the statement
of Theorem \ref{theorem:BCRW}. Then for $y_1<\dots<y_{2n}$
 \[
 \Phfac(y)   \mathcal{L} \left( \Sigma_y\right) (\eta) 
 = \sum_{i=1}^{2n} (L^{p \phi, q \phi}_{y_i} - c_0) \left( \Phfac(y) \Sigma_y(\eta) \right)
\]
where $c_0 = \frac12 (r+l) - (p+q)(\phi-1) = \frac{p+q}{2} \, (\phi-1)^2 \geq 0 $.
\end{Lemma}

\textbf{Proof.} This is a direct calculation. For a function of one variable we find that a change
of $f: \Z \to \R$ to $\tilde{f}(y) = c^y f(y)$, for $c \neq 0$,  produces the change
\[
L^{p,q} \tilde{f}(y) = c^y L^{pc^{-1},qc}f(y)  + (pc^{-1} +qc -p-q) c^y f(y).
\]
We apply this in the even coordinates with $c=\phi$ and in the odd co-ordinates with $c= \phi^{-1}$.
The value of $\phi$ is chosen so that the corresponding difference operators now agree on both 
sets of coordinates. Different potential terms are created at odd or even co-ordinates, but these can be 
summed and then 
redistributed equally between all co-ordinates, which yields the constant $c_0$.  
The equivalent expressions for $c_0$ follow from the definition of $\phi$. \qed
\begin{Lemma} \label{lemma:emptyintervalpfaffian}
 For all $\eta\in\{0,1\}^\Z$, for all $n\ge1$, $y_1\le\dots\le y_{2n}$ and $t\ge0$
 \[
\Phfac(y)   \E_\eta\left[\Sigma_y(\eta_t)\right] = \Pf (K^{(2n)}(t,y)),
 \]
 where $K^{(2n)}(t,y)$ is the anti-symmetric $2n\times2n$ matrix with entries $K_t(y_i,y_j)$ defined, for $i<j$,  by 
\[
K_t(y,z) = \phi^{z-y} \, \Pp_{\eta}[ \eta_t[y,z)=0], \quad \mbox{for $t \geq 0$ and $y \leq z$.}
\]
\end{Lemma} 

\textbf{Proof of Lemma~\ref{lemma:emptyintervalpfaffian}.}
We follow closely the arguments for the pure coalescing case, as in Lemma 7 of \cite{GPTZ:15}, pointing out the changes needed here but 
refer the reader to \cite{GPTZ:15} details. We recall the notation for the discrete cells $V_{2n}$ and their faces 
$(\partial V^{(i)}_{2n}: 1 \leq i \leq 2n-1)$, defined as follows:
\begin{eqnarray*}
V_{2n} & = & \{ y \in \Z^{2n}: y_1<y_2<\ldots < y_{2n}\}, \\
\partial V^{(i)}_{2n} & = & \{ y \in \Z^{2n}: y_1<y_2<\ldots < y_i = y_{i+1} < \ldots < y_{2n}\}.
\end{eqnarray*}
We also use the notation $y^{i,i+1}$ for the vector $y$ with coordinates $y_i$ and $y_{i+1}$ removed; thus 
when $n \geq 2$, for  $y \in \partial V^{(i)}_{2n}$ we have  $y^{i,i+1} \in V_{2n-2}$ 

To establish the identity in the Lemma, one checks that  both sides are solutions to the system of ODEs, in this case
\[
(\mbox{ODE})_{2n} \hspace{.3in} \left\{ \begin{array}{rcll}
\partial_t u^{(2n)}(t,y) &=& \sum_{i=1}^{2n} \, (L^{p\phi,q\phi}_{y_i} - c_0) u^{(2n)}(t,y) & \mbox{on $[0,\infty) \times V_{2n}$,} \\
u^{(2n)}(t,y) & = & u^{(2n-2)}(t,y^{i,i+1}) &  \mbox{on $[0,\infty) \times \partial V_{2n}^{(i)}$,} \\
u^{(2n)}(0,y) & = & \Phi(y) \Sigma_y(\eta) &  \mbox{on $V_{2n}$,} 
\end{array} \right.
\] 
taking $u^{(0)} \equiv 1$.
This infinite system can be shown by induction on $n$ to have unique solutions, within the class
of functions that suitable exponential growth at infinity. As in \cite{GPTZ:15}, the fact that $(t,y) \mapsto\E_\eta\left[ \Phfac(y) \Sigma_y(\eta_t)\right]$ is a solution follows from Lemma \ref{lemma:balancing}
and the extra fact that the phase factor satisfies $\Phi(y) = \Phi(y^{i,i+1})$ on the boundary $\partial V_{2n}^{(i)}$.

The fact that the Pfaffian is also the solution to this system follows as in the non-branching case in \cite{GPTZ:15}. 
The starting point is that the function $K_t(y,z)$ solves the equation
\begin{equation} \label{K2new}
\partial_t K(y,z) = \left( L^{p\phi,q\phi}_y +  L^{p\phi,q\phi}_z - 2 c_0 \right) K_t(y,z) \quad \mbox{for $(y,z) \in V_2$}
\end{equation}
with the boundary condition $K_t(y,y) = 1$ for all $t >0$. The Pfaffian is made up terms that are products of copies of $K$
using all the variables $y_1,\ldots,y_{2n}$, each of which can be checked to be a solution of the the full system. The combination
of terms in the Pfaffian is exactly what is needed t satify the boundary conditions of the system.
The details of these arguments are explained in Lemma 7 of \cite{GPTZ:15} with the only change being that we need to verify the 
extra phase term does not affect the  initial condition being satisfied.  However we may rewrite the entries in the 
Pfaffian at time zero using
\[
K_0(y,z) = \phi^{z-y} \, \ind(\eta[y,z)=0) = 
\lim_{\theta\downarrow0}\frac{\theta^{\eta[a,z)}}{\theta^{\eta[a,y)}}\frac{\phfac^{z-a}}{\phfac^{y-a}} \quad \mbox{for $a<y<z$.}
\]
The Pfaffian $\Pf (K^{(2n)}(0,y))$ is therefore the limit of Pfaffians of a matrix $A$ with entries in quotient form  $A_{ij} = a_i/a_j$ for $i<j$.
For such matrices $\Pf(A) = a_2 a_4 \ldots a_{2n}/a_1 a_3 \ldots a_{2n-1}$ (see the appendix of \cite{GPTZ:15} for example) and hence, taking
$a< \min\{y_i\}$, 
\[
 \Pf(K^{(2n)}(0,y))
 =\lim_{\theta\downarrow0}\prod_{i=1}^n\frac{\theta^{\eta[a,y_{2i})}}{\theta^{\eta[a,y_{2i-1})}}\frac{\phi^{y_{2i}-a}}{\phi^{y_{2i-1}-a}}
=\Phi(y)  \Sigma_y(\eta),
\]
as required. \qed

\textbf{Remark. 2.1.} The last lemma is the point at which to observe that for certain random initial conditions, the 
Pfaffian property is still true. Indeed suppose that $\eta_0$ is random but that 
$\E \left[\Phfac(y) \Sigma_y(\eta_0)\right] $ is still given by a $2n \times 2n$ Pfaffian with entries
$K_0(y_i,y_j)$ for $i<j$, for some $K_0$ of exponential growth. The statement of the lemma then still holds,
and so does Theorem \ref{theorem:BCRW}, which is deduced from the lemma without any changes.
A simple example is when $\eta_0(x)$ are independent Bernoulli($\theta_x$) variables. Then 
the condition above is true with
\[
K_0(y,z) = \phi^{z-y} \prod_{k \in [y,z)} (1-\theta_k).
\]
A similar observation holds for the AWRPI model discussed in the next section. 

\textbf{Proof of Theorem~\ref{theorem:BCRW}.}
As in \cite{GPTZ:15}, the desired particle intensities $\E_\eta\left[\eta_t(x_1)\dots\eta_t(x_n)\right]$ 
may be recovered from the empty interval probabilities. 
% The formula
% \[
%  D^+_z\sigma_{y,z}(\eta) = -\eta(z)\sigma_{y,z}(\eta)
% \]
% implies that
% \[
%  D^+_z\sigma_{y,z}(\eta)|_{z=y} = -\eta(z).
% \]
From Lemma~\ref{lemma:emptyintervalpfaffian} we have
%\begin{equation} \label{eq:emptyintformula}
\[
 \E_\eta \left[\Sigma_y(\Eta_t)\right] =  \Phi(y)^{-1} \Pf(K^{(2n)}(t,y)).
\]
%\end{equation}
%Remark on thinning interpretation of this formula? This leads to immediate proof of various 
%conjugation steps below...but do we know there is a point process before thinning? 
The factor $\Phi^{-1}(y)$ can be expressed as the determinant of a diagonal matrix $D(y)$ with entries $D_{ii}(y) = \phi^{(-1)^{i+1} y_i}$ for $i=1,\ldots,2n$.
The empty interval probabilities can then be expressed as a single Pfaffian
\begin{equation} \label{eq:emptyintformula}
 \E_\eta \left[\Sigma_y(\Eta_t)\right] =   \Pf( D(y) K^{(2n)}(t,y) D(y)).
\end{equation}
Note the $ij$'th entry of the matrix $D(y) K^{(2n)}(t,y) D(y)$ is still a function only of the variables $y_i,y_j$. 
We now follow the argument in \cite{GPTZ:15}, where the intensities are derived from the empty interval probabilities via 
discrete derivatives. This leads to the process $\eta_t$  being a Pfaffian point 
process with a kernel $\hat{\K}(y,z)$ where for $y<z$
\[
  \left(\begin{array}{cc}
         \phi^{y+z}  K_t(y,z) & - D^+_{z}\left(\phi^{y-z} K_t(y,z)\right) \\
         - D^+_{y}\left(\phi^{z-y} K_t(y,z)\right) & D^+_{y}D^+_{z}\left(\phi^{-y-z} K_t(y,z)\right)
        \end{array}\right),
\]
and 
\[
\hat{\K}_{12}(y,y) = -D^+_z \left(\phi^{y-z} K_t(y,z)\right)|_{z=y} = 1- \phi^{-1}  K_t(y,y+1).
\] 
It remains to massage this kernel $\hat{\K}$ into the form $\K$ stated in the Theorem, which uses only row and column operations that
can be represented by conjugation with suitable matrices, that is we may replace $\hat{\K}(y,z)$ by $A(y)\hat{\K}(y,z) A^T(z)$
for any $2$-by-$2$ matrix $A(y)$ (depending only on the variable $y$) that has determinant $1$.  
% One aim of this manipulation is to show that the final kernel is still in derived form 

Expanding out the discrete derivatives in $\hat{\K}$ using the discrete product rule, and then conjugating the final matrix with a
diagonal matrix 
$A(y) = \left( \begin{array}{cc} \phi^{-y} & 0 \\ 0 & \phi^y \end{array} \right)$  leads to an equivalent kernel, which we still denote $\hat{\K}$, with entries
\[
 \begin{array}{rcl}
   \hat{\K}_{11}(y,z) & = & K_t(y,z); \\
   \hat{\K}_{12}(y,z) & = & -\left( \phi^{-1} K_t(y,z+1) - K_t(y,z) \right);\\
   \hat{\K}_{21}(y,z) & = & -\left( \phi^{-1} K_t(y+1,z) - K_t(y,z) \right);\\
   \hat{\K}_{22}(y,z) & = & \phi^{-2} K_t(y+1,z+1) - \phi^{-1} K_t(y,z+1) - \phi^{-1} K_t(y+1,z) + K_t(y,z),
  \end{array}
\]
$\hat{\K}_{12}(y,y)=1- \phi^{-1} K_t(y,y+1)$. 

Subtracting the first row and column from the second row and column, and then further conjugating with a diagonal matrix
$A(y) = \left( \begin{array}{cc} \phi^{-1/2} & 0 \\ 0 & \phi^{1/2} \end{array} \right)$
 gives the equivalent kernel $\K$
\[
\hat{\K}(y,z) = \phi^{-1}
\left(\begin{array}{cc}
         K_t(y,z) & - K_t(y,z+1) \\
         -K_t(y+1,z) & K_t(y+1,z+1)
        \end{array}\right),
\]
with $\hat{\K}_{12}(y,y)=1- \phi^{-1} K_t(y,y+1)$. 
Finally, the desired kernel $\K$ is obtained by again subtracting the first row and column from the second. \qed

\textbf{Remarks.} \textbf{2.2.} Letting $t \to \infty$ the process, for any non-zero initial condition, converges to
an invariant Bernoulli distribution. It is fun to see this via the Pfaffian kernels by
solving for $K_t(y,z)$ explicitly.
Consider first the maximal initial condition $\eta_0(x) =1$ for all $x$. Then $K_0(y,z) = 0$
and the kernel $K_t(y,z)$ has a representation in terms of a pair of independent continuous time random walkers $(Y_t,Z_t)$ with 
generator $q \phi D^+ + p \phi  D^-$, started at $Y_0=y, Z_0=z$. Let $\tau = \inf\{t:Y_t=Z_t\}$. Then the 
equation (\ref{K2new}) solved by $K_t(y,z)$ implies that 
$M_s := K_{t-s}(Y_s,Z_s) \exp(-2 c_0 s)$ is a martingale for $s \in [0,t \wedge \tau]$
and hence that
\[
K_t(y,z) = M_0 = \E[M_{t \wedge \tau}] = \E \left[ e^{-2c_0 \tau} \ind(\tau \leq t)\right]] \uparrow K_{\infty}(y,z) = \E[ e^{-2 c_0 \tau}] .
\]
Solving explicitly we find $K_{\infty}(y,z) = \phi^{-(z-y)}$, and the Pfaffian kernel of the form 
(\ref{discretekernel1}) corresponding to $K_{\infty}(y,z)$  is 
\[
  \K(y,z) = \frac{\phi^{y-z-1}}{2}\left(\begin{array}{cc}
	      1 & (1-\phi^{-1}) \\
	    (1-\phi) & (1-\phi)(1-\phi^{-1})
	    \end{array}\right),
\]
and $\K_{12}(y,y) = 1-\phi^{-2}$.
A little manipulation (using row and column operations as in (\ref{star1}) for the ARWPI model kernel) shows that this is a kernel for a product Bernoulli($1- \phi^{-2}$) distribution. Convergence of $K_t(y,z)$ implies that all entries in the Pfaffian kernel converge,
which in turn implies that the process $\eta_t$ converges as $ t \to \infty$ to the product Bernoulli($1-\phi^{-2}$) in distribution in the product topology.

For general non-zero initial conditions the same is true. Rather than analyse the kernel, we use a simple coupling argument for attractive nearest neighbour systems. All non zero solutions can be coupled between the maximal solution 
and a solution started from a single point. It therefore is enough to prove convergence for the process $\eta^0_t$ started from a single occupied site, say the origin. But this process can be coupled with the process $\eta^{\Z}_t$ started from all occupied sites.
Indeed by a graphical construction (or equivalently solving a system of equations using the same Poisson drivers for jump and branch events) shows that 
\[
\eta^0_t(y) = \eta^{\Z}_t(y) \quad \mbox{for all $y \in [l_t,r_t]$,}
\]
where $l_t,r_t$ mark the leftmost and rightmost occupied site in $\eta^0_t$. The behaviour of the pair $\{l_t,r_t\}$ is however easy to understand: provided $p+l>q$ and $q+r>p$ we can ensure $l_t \to -\infty$ and $r_t \to \infty$. Under these conditions
the process looks like $\eta^Z_t$ in a growing interval, and we already know $\eta^{\Z}_t$  converges to Bernoulli equilibrium. 

\textbf{2.3.} It is natural to look for a spatially inhomogeneous version of the BCRW model, where
$p_x,q_x,l_x,r_x$ are allowed to be site dependent. This was explored in the thesis \cite{thesis} and 
the Pfaffian property can survive, but under a somewhat stronger condition on the parameters that we do not 
yet fully understand.
\section{Proof of Theorem \ref{theorem:ARWPI}.}  \label{s3}
The result for the annihilating model with immigration follows by very similar lines, and we remark only on the changes
caused by the new immigration term. The result holds for systems with spatially inhomogeneous rates. 
There is no additional complexity in the proof, so we continue in this general framework.

The generator for (spatially inhomogeneous) ARWPI is given, for suitable $F:\{0,1\}^\Z\to\R$ by
\begin{eqnarray} 
 \mathcal{L}  F(\eta) & = & \sum_{x\in\Z} q_x\left (F(\eta_{x,x-1})-F(\eta)  \right) + \sum_{x\in\Z} p_x\left(F(\eta_{x-1,x}) - F(\eta)\right) \nonumber \\
&& \hspace{.2in} + \sum_{x\in\Z} m_x\left (F(\eta^{i}_{x-1,x})-F(\eta) \right), 
\label{eq:processgeneratorarwpi}
\end{eqnarray}
where $\eta_{x,y}$ is the new configuration resulting from a jump from $x$ onto $y$, that is 
\[
\left\{
\begin{array}{l}
\eta_{x,y}(z) = \eta(z) \quad \mbox{for $z \not \in \{x,y\}$,} \\ 
\eta_{x,y}(x) = 0, \quad  
\eta_{x,y}(y) = \eta(x) + \eta(y) \mod(2),
\end{array} \right.
\]
and where $\eta^{i}_{x-1,x}$ is the new configuration resulting from a pair immigration onto $\{x-1,x\}$ defined by
\[
\left\{
\begin{array}{l}
\eta^{i}_{x-1,x}(z) = \eta(z) \quad \mbox{for $z \not \in \{x-1,x\}$,} \\ 
\eta^i_{x-1,x}(z) = 1-\eta(z) \quad \mbox{for $z \in \{x-1,x\}$}.
\end{array} \right.
\]
Note that any immigrating particle instantly annihilates with any existing particle on the target site. We suppose
$m_x,p_x,q_x$ are uniformly bounded, so that this generator uniquely determines a Markov process. 

The following spin product function is a Markov duality function for this generator (as used for the pure annihilating model in 
\cite{GPTZ:15}). For $n\ge1$ and $y=(y_1,\dots,y_{2n})$ with $y_1\le\dots\le y_{2n}$ we define the spin product
\[
 \Sigma_y(\eta) = \prod_{i=1}^n(-1)^{\eta[y_{2i-1},y_{2i})}.
\]
We define the one-particle operator $L$, acting on $f:\Z\to\R$, by
\begin{equation}\label{eq:Loperator}
  Lf(x) = q_xD^+f(x) + p_xD^-f(x) - 2m_xf(x).
\end{equation}

\begin{Lemma}\label{lemma:MarkovdualityA}
 For $y_1<\dots<y_{2n}$ the action of the generator $\mathcal{L}$ on $\Sigma_y(\eta)$ is given by
 \[
  \mathcal{L}\Sigma_y(\eta) = \sum_{i=1}^{2n}L_{y_i}\Sigma_y(\eta).
 \]
\end{Lemma}

\textbf{Proof of Lemma~\ref{lemma:MarkovdualityA}.}
As in \cite{GPTZ:15} the terms of $\mathcal{L}$ coming from particle motion contribute 
\[
  \sum_{x\in\Z} q_x \left (\Sigma_y(\eta_{x,x-1})-\Sigma_y(\eta)  \right) +\sum_{x\in\Z} p_{x} \left(\Sigma_y(\eta_{x-1,x}) - \Sigma_y(\eta)\right).
\]
It remains to compute the immigration term. Note that the modified immigration configuration $\eta^i_{x-1,x}$ differs from $\eta$ on at most two sites, $x-1$ and $x$. Since the $y_i$ are strictly ordered, the intervals $[y_{2i-1},y_{2i})$ are separated by at least one site, whereby a pair of adjacent sites $-1,x$ can intersect at most one of the intervals. In particular any change due to immigration affects at most one interval $[y_{2i-1},y_{2i})$, leading to the formula
\[
 \Sigma_y(\eta^{i}_{x-1,x})-\Sigma_y(\eta) 
 = \sum_{i=1}^n \left(\prod_{j\neq i} (-1)^{\eta[y_{2j-1},y_{2j})}\right) 
\left((-1)^{\eta^i_{x-1,x}[y_{2i-1},y_{2i})} - (-1)^{\eta[y_{2i-1},y_{2i})} \right).
\]
Fix $y<z$ and consider the generator contribution for a single spin product $(-1)^{\eta[y,z)}$, namely
\[
 \sum_{x\in\Z} m_x \left( (-1)^{\eta^{i}_{x-1,x}[y,z)} - (-1)^{\eta[y,z)} \right).
\]
The terms indexed by $x\le y-1$ and $x\ge z+1$ are zero, as the modified configuration is unchanged in the interval $[y,z)$. The terms $y+1\le x\le z-1$ are also zero, since
the immigration of two particles does not change the parity of $\eta[y,z)$. 
The remaining terms give identical non-zero contributions: for $x=y$ or $x=z$ 
\[
(-1)^{\eta^{i}_{x-1,x}[y,z)} - (-1)^{\eta[y,z)}
 = \prod_{\substack{w=y\\w\neq x}}^{z-1}(-1)^{\eta(w)} \left((-1)^{1-\eta(x)} - (-1)^{\eta(x)}\right)
 = -2 (-1)^{\eta[y,z)}.
\]
All together the immigration term is given by
\[
 \sum_{x\in\Z} m_x \left( \Sigma_{y}(\eta^{i}_{x-1,x}) - \Sigma_{y}(\eta) \right)
 = -2\Sigma_{y}(\eta)\sum_{i=1}^{2n}m_{y_i}.
\]
Collecting the jump and immigration terms gives the desired expression. \qed

\textbf{Proof of Theorem~\ref{theorem:ARWPI}.}
Following the argument from \cite{GPTZ:15}, we first claim that for all $\eta\in\{0,1\}^\Z$, for all $n\ge1$, $y_1\le\dots\le y_{2n}$ and $t\ge0$\[
  \E_\eta\left[\Sigma_y(\eta_t)\right] = \Pf (K^{(2n)}(t,y)),
 \]
 where $K^{(2n)}(t,y)$ is the anti-symmetric $2n\times2n$ matrix with entries $K_t(y_i,y_j)$ for $i<j$, defined by 
$K_t(y,z) = \E_{\eta}[ (-1)^{\eta_t[y,z)}]$. The particle intensities $\E_{\eta} \left[\eta_t(x_1)\dots\eta_t(x_n)\right]$ can then be recovered from product spin expectations via discrete derivatives and yield the stated kernel $\K(y,z)$.  \qed
\section{Some continuum Pfaffian point process limits. }  \label{s4}
The entries for the Pfaffian kernels $\K(x,y)$ in both the branching model and the immigration model, are determined by a scalar function
$K_t(x,y)$ that solves a certain discrete heat equation. Under diffusive space time scaling, and with suitable scaling of the parameters, we can obtain natural 
limiting Pfaffian point processes $X^{(c)}$ on $\R$, with associated continuum kernels $\K^{(c)}$ (where the superscript (c) stands for continuum). 
We record here certain examples, simply to add to the rather small current list of explicit Pfaffian point process kernels. Two points are perhaps of greatest interest:

\textbf{1.}  Unlike the continuum examples from \cite{GPTZ:15}, alongside the diffusive scaling of space-time, the reaction parameters controlling branching and immigration must be simultaneously scaled, so that they have a non-trivial effect on the continuum limit. Indeed branching but instantly coalescing Brownian motions do not have a simple meaning, and nor does
immigration of instantly annihilating pair of Brownian motions onto the same point. 
However, since both discrete process are Pfaffian whose entire statistics are controlled
by a kernel whose entries solve a discrete PDE, the correct scaling for the parameters is easily revealed by examining the convergence for the differential equations.

\textbf{2.} The discrete equations behind coalescing models with branching and annihilating models with pairwise immigration are both discrete heat equations with a constant potential. This can be used to show there is an equality in law for the fixed time particle positions between these two models, if parameter values and initial values are chosen carefully. This connection exists at the discrete level (see \cite{thesis}) but is most transparent for the limiting continuum models, and we detail this in the remarks after example (b). 

%\textbf{Examples of continuum kernels.} 

In each of the four examples below we define
\[
X^{(\epsilon)}_t(dx) = \eta_{\epsilon^{-2}t}(\epsilon^{-1}dx) \quad \mbox{on $\epsilon \Z$}
\]
where $\eta_t$ is one of the models studied earlier, with an initial condition and $\epsilon$ dependent parameters which we will specify.
The point process $X^{(\epsilon)}_t$ will be a Pfaffian point process on $\epsilon \Z$ with a kernel $\K^{(\epsilon)}_t$. The diffusive scaling above is 
chosen so that an isolated non-interacting particle will converge to a Brownian motion. 
We claim convergence of the particle system only at a fixed time. Indeed, for all $t \geq 0$, in each 
of the examples (a),(b),(c) below
we claim $X_t^{(\epsilon)} \to X^{(c)}_t$ in distribution, on the space of locally finite point measures on $\R$ with 
the topology of vague convergence (for the final example (d) we restrict to a region away from the origin). 
Moreover the limit $X^{(c)}_t$ is a simple point process and a Pfaffian point process on $\R$. In our examples we can often solve explicitly for the limiting kernel $\K^{(c)}_t(x,y)$.

The entries of $\K^{(\epsilon)}_t$ are in terms of a scalar function
$K^{(\epsilon)}_t(y,z)$ that will solve a lattice PDE that naturally scales to a continuum PDE. 
The proof of the convergence $X^{(\epsilon)}_t \to X_t$ follows from the suitable convergence of 
$K^{(\epsilon)}_t(y,z)$ and their discrete derivatives to the analogous solutions of a continuum PDE, by plugging in to the
kernel convergence Lemma 9 from \cite{GPTZ:15}. However we omit the details verifying all the conditions of this lemma.
% More here?
%

\textbf{(a) Annihilating model with constant pairwise immigration.}

We consider the ARWPI model with parameters $p_x=q_x=\alpha>0$ and $m_x = \beta \epsilon^{-2} \geq 0$ for all $x$, and with zero initial condition. 
From Theorem \ref{theorem:ARWPI} the process $X^{(\epsilon)}_t$ is a Pfaffian point process on $\epsilon \Z$ with kernel 
$\K_t^{(\epsilon)}$ of the form
 \begin{eqnarray} 
&& \K^{(\epsilon)}_t(y,z)
 = \frac{\epsilon}{2} \left(\begin{array}{cc}
           K^{(\epsilon)}_t(y,z) & -D^{(\epsilon)}_2 K^{(\epsilon)}_t(y,z) \\
           - D^{(\epsilon)}_1 K^{(\epsilon)}_t(y,z) & D^{(\epsilon)}_1 D^{(\epsilon)}_2 K^{(\epsilon)}_t(y,z)
           \end{array}\right)
           \quad \mbox{for $y<z$,}    \nonumber \\
 && (\K^{(\epsilon)}_t)_{12}(y,y) = - \frac{\epsilon}{2} D^{(\epsilon)}_2 K^{(\epsilon)}_t(y,y)
 \label{epskernela}
 \end{eqnarray}
where $D^{(\epsilon)}_i$ is the right discrete derivative on $\epsilon \Z$ (that is
$D^{(\epsilon)} f(x) = \epsilon^{-1}(f(x+\epsilon) - f(x))$) acting on the $i$'th variable.
The function $K^{(\epsilon)}_t(y,z)$ solves, for $y,z \in \epsilon \Z$ with $y \leq z$, 
\begin{equation}
\left\{ \begin{array}{rcll}
\partial_t K^{(\epsilon)}_t &=&  \alpha \Delta^{(\epsilon)} K^{(\epsilon)}_t - 2\beta K^{(\epsilon)}_t,\\
K^{(\epsilon)}_t(y,y) & = & 1, \\
K^{(\epsilon)}_0(y,z) & = & 1.
\end{array} \right. \label{epspdea}
\end{equation}
Here $\Delta^{(\epsilon)}$ is the discrete Laplacian on $(\epsilon \Z)^2$.

The limit $X_t$ is Pfaffian on $\R$ with kernel of the form 
 \begin{eqnarray} 
&& \K_t^{(c)}(y,z)
  = \frac12 \left(\begin{array}{cc}
           K^{(c)}_t(y,z) & -D_2 K^{(c)}_t(y,z) \\
           - D_1 K^{(c)}_t(y,z) & D_1 D_2 K^{(c)}_t(y,z)
           \end{array}\right)
 \quad \mbox{for $y<z$,} \nonumber \\
 && (\K_t^{(c)})_{12}(y,y) = - \frac12 D_2 K^{(c)}_t(y,y) \label{ctskernela}
 \end{eqnarray}
 ($D_i$ is the derivative in the $i$th co-ordinate)
where $K^{(c)}_t(y,z)$ is $C^2$ on $\{y,z\in \R^2:y \leq z\}$ and solves
\begin{equation}
\left\{ \begin{array}{rcl}
\partial_t K^{(c)}_t(y,z) &=&  \alpha \Delta  K^{(c)}_t(y,z)  - 2\beta K^{(c)}_t(y,z)\\
K^{(c)}_t(y,y) & = & 1 \\
K^{(c)}_0(y,z) & = & 1.
\end{array} \right. \label{ctspdea}
\end{equation}
The unique bounded solution $K^{(c)}_t(y,z)$ has a probabilistic representation in terms of a two dimensional Brownian motion $(Y_t,Z_t)$, run at 
rate $2 \alpha$ (that is scaled to have variance  $2 \alpha t$ at time $t$) and started at $(y,z)$, 
namely
\[
K^{(c)}_t(y,z) = \E \left[ e^{-2\beta (t \wedge \tau)} \right] \quad \mbox{where $\tau = \inf\{t:Y_t=Z_t\}$.} 
\]
Solving for $K^{(c)}_t(y,z)$ explicitly allows one to read off the one point density 
\[
\rho^{(1)}_t(y)= - \frac12 D_2 K^{(c)}_t(y,y) = \frac12 \sqrt{\frac{\beta}{\alpha}} 
\erf(\sqrt{2t\gamma})
\]
(where the error function is defined by $\erf(x) = (2/\sqrt{\pi}) \int^x_0 \exp(-t^2) dt$.) 
The kernel also has a limit as $t \to \infty$, in particular $K^{(c)}_t(y,z) \to K^{(c)}_{\infty}(y,z)$ where
\[
K_{\infty}^{(c)}(y,z) = \E[ e^{-2\beta \tau}] = e^{-\sqrt{\frac{\beta}{\alpha}}(z-y)}.
\]
It is no longer enough in the continuum to just examine convergence of $K^{(c)}_t$, but an examination of the 
exact formula shows that both $K^{(c)}$ and its first two derivatives converge, uniformly over $y,z$, as $t \to \infty$, which implies that the continuum point processes $X^{(c)}_t$ converge as $t \to \infty$ to a point process $X^{(c)}_{\infty}$ 
(one can follow the steps of the proof of Lemma 4 from \cite{GPTZ:15}). The limit has 
kernel
 \begin{equation}
 \K^{(c)}_{\infty}(y,z) = 
 \frac{1}{2}\left(\begin{array}{cc}
   e^{-\sqrt{\frac{\beta}{\alpha}}(z-y)} & \sqrt{\frac{\beta}{\alpha}}
   e^{-\sqrt{\frac{\beta}{\alpha}}(z-y)} \\
    -\sqrt{\frac{\beta}{\alpha}}e^{-\sqrt{\frac{\beta}{\alpha}}(z-y)} & -\frac{\beta}{\alpha}e^{-\sqrt{\frac{\beta}{\alpha}}(z-y)}
 \end{array}\right) \quad \mbox{for $y<z$,} \label{poissonkernel}
 \end{equation}
 and $\K^{(c)}_{\infty,12}(y,y)=\frac{1}{2}\sqrt{\frac{\beta}{\alpha}}$. 
One can identify the limit, this time a disguised form for the kernel for a Poisson process. Indeed the same row and column operations as in the discrete case allow the Pfaffian of the above kernel to be easily computed explicitly and the $n$-point intensity is given by
\[
 \rho^{(n)}_\infty(y_1,\dots,y_n) \equiv \left(\frac{1}{2}\sqrt\frac{\beta}{a}\right)^n.
\]
Thus the distribution of the (continuum) point process in the large time limit is a 
Poisson process rate $\frac{1}{2}\sqrt{\frac{\beta}{\alpha}}$. 
The four entries in the kernel converge exponentially to the $t=\infty$ limit. 
As in the discrete ARWPI model, this can be used to show the exponentially fast convergence $\rho^{(n)}_t(y_1,\ldots,y_n) \to \rho^{(n)}_{\infty}(y_1,\ldots,y_n)$ as $t \to \infty$, for any fixed $n$ and uniformly over $y_i$.  

Explicit formulae can be found for a variety of other initial 
conditions (see \cite{thesis}). For example for an initial Bernouilli($ \epsilon \theta$) condition, where $\theta$ is fixed,  only the initial condition in (\ref{epspdea}) changes to $K^{(\epsilon)}_0(y,z) = (1-2 \epsilon \theta)^{\epsilon^{-1}(z-y)}$, and the 
initial condition for the limiting PDE (\ref{ctspdea}) changes to $K^{(c)}_0(y,z) = \exp(-2 \theta (z-y))$. In the maximal case,
$\eta_x=1$ for all $x$, the initial conditions $K^{(\epsilon)}_0(y,z) = (-1)^{\epsilon^{-1}(z-y)}$ are extermely oscilatory, but they 
converge in distribution to the zero function which is sufficient to imply suitable convergence of the kernels at a fixed times $t>0$.  

\textbf{(b) Branching and coalescing model with maximal initial condition.} 

We consider the BCRW model with parameters $p_x=q_x=\alpha>0$ and $l = r = 2 \epsilon \sqrt{\alpha \beta}>0$, 
and with maximal initial condition, that is $\eta(x) =1$ for all $x$. (We have chosen the form of the branching rate parameters so that the limit has a convenient form). 
From Theorem \ref{theorem:ARWPI} the process $X^{(\epsilon)}_t$ is a Pfaffian point process on $\epsilon \Z$ with kernel 
$\K_t^{(\epsilon)}$ of the form
 \begin{eqnarray} 
&& \K^{(\epsilon)}_t(y,z)
 = \epsilon \phi \left(\begin{array}{cc}
           K^{(\epsilon)}_t(y,z) & -D^{(\epsilon)}_2 K^{(\epsilon)}_t(y,z) \\
           - D^{(\epsilon)}_1 K^{(\epsilon)}_t(y,z) & D^{(\epsilon)}_1 D^{(\epsilon)}_2 K^{(\epsilon)}_t(y,z)
           \end{array}\right)
           \quad \mbox{for $y<z$,}    \nonumber \\
 && (\K^{(\epsilon)}_t)_{12}(y,y) = 1 - \phi^{-1} K^{(\epsilon)}_t(y,y+\epsilon).
 \label{epskernelb}
 \end{eqnarray}
The function $K^{(\epsilon)}_t(y,z)$ solves, for $y,z \in \epsilon \Z$ with $y \leq z$, 
\begin{equation}
\left\{ \begin{array}{rcll}
\partial_t K^{(\epsilon)}_t &=&  \alpha \phi \Delta^{(\epsilon)} K^{(\epsilon)}_t - 2 \epsilon^{-2} c_0 K^{(\epsilon)}_t,\\
K^{(\epsilon)}_t(y,y) & = & 1, \\
K^{(\epsilon)}_0(y,z) & = & 0.
\end{array} \right. \label{epspdeb}
\end{equation}
Examination of the constants $\phi$ and $c_0$ shows that
\[
\phi = 1 + \epsilon \sqrt{\beta/\alpha} - \epsilon^2 (\beta/2 \alpha) + O(\epsilon^3), \qquad
c_0 = \beta \epsilon^2 + O(\epsilon^3). 
\]
The limit $X_t$ is Pfaffian on $\R$ with kernel of the form 
 \begin{eqnarray} 
&& \K_t^{(c)}(y,z)
  = \left(\begin{array}{cc}
           K^{(c)}_t(y,z) & -D_2 K^{(c)}_t(y,z) \\
           - D_1 K^{(c)}_t(y,z) & D_1 D_2 K^{(c)}_t(y,z)
           \end{array}\right)
 \quad \mbox{for $y<z$,} \nonumber \\
 && (\K_t^{(c)})_{12}(y,y) = -  D_2 K^{(c)}_t(y,y) + \sqrt{\beta/\alpha} \label{ctskernelb}
 \end{eqnarray}
where $K^{(c)}_t(y,z)$ is $C^2$ on $\{y,z\in \R^2:y \leq z\}$ and solves
\begin{equation}
\left\{ \begin{array}{rcl}
\partial_t K^{(c)}_t(y,z) &=&  \alpha \Delta  K^{(c)}_t(y,z)  - 2\beta K^{(c)}_t(y,z)\\
K^{(c)}_t(y,y) & = & 1 \\
K^{(c)}_0(y,z) & = & 0.
\end{array} \right. \label{ctspdeb}
\end{equation}
Note that the term $(\K_t^{(c)})_{12}(y,y)$ requires a bit more care than in example (a) and that
an extra term $\sqrt{\beta/\alpha} $ emerges, as follows: 
\begin{eqnarray*}
\epsilon^{-1} K^{(\epsilon)}_t(y,y) & = & \epsilon^{-1} \left( 1 - \phi^{-1} K^{(\epsilon)}_t(y,y+ \epsilon) \right) \\
& = &  \epsilon^{-1} \left( 1 - (1  + \epsilon \sqrt{\beta/\alpha} + O(\epsilon^2))^{-1} 
(1 + \epsilon D^{(\epsilon)}_2 K^{(\epsilon)}_t(y,y)) \right) \\
& = & - D^{(\epsilon)}_2 K^{(\epsilon)}_t(y,y) + \sqrt{\beta/\alpha} + O(\epsilon^2) \\
& \to & - D_2 K^{(c)}_t(y,y) + \sqrt{\beta/\alpha} \quad \mbox{as $\epsilon \downarrow 0$.}
\end{eqnarray*}

The $t \to \infty$ limit follows the lines of the previous example, and $X^{(c)}_t \to X^{(c)}_{\infty}$ 
where the limit has Pfaffian kernel that is twice the one in (\ref{poissonkernel}), and with the extra difference that
\[
(\K^{(c)}_{\infty})_{12}(y,y)= -D_2 K_{\infty}^{(c)}(y,y) + \sqrt{\frac{\beta}{\alpha}} = 2 \sqrt{\frac{\beta}{\alpha}}.
\]
Similar row and column manipulations as in the discrete case show that this kernel encodes a Poisson process of rate $2 \sqrt{\frac{\beta}{\alpha}}$.

\textbf{Remarks.} \textbf{4.1.} 
In many formulations of determinantal point processes, the determinantal kernel $D$ is associated to 
an integral operator $D$ on $L^2(\R)$, and the diagonal values $D(y,y)$ are linked to those of
$(D(y,z):y<z)$ by the fact that the operator is assumed to be of trace class. One might ask the same 
for the Pfaffian case, asking for four operators $\K_{ij}$ on $L^2(\R)$. 
For our examples this link is broken: the operators $\K_{ij}$ acting on $L^2(\R)$ would have 
discontinuities along $y=z$ and are not expected to be trace class. The diagonal values 
$\K_{12}(y,y)$ are not given, for example, as even the continuous limit of $\K_{12}(y,z)$.
This is also the case for classical Pfaffian kernels, for example for GOE.

The operator formulation is useful, for example when applying the theory
of Fredholm determinants or Fredholm Pfaffians, and in classification theorems.  
However, we state our continuum Pfaffian kernels in the form of the 
five measurable functions, namely $(\K_{ij}(y,z): y,z \in \R, y<z))_{i.j \in \{1,2\}}$ and $(\K_{12}(y,y): y \in \R)$. 
These five functions are what appear in the Pfaffian formulae for the intensities $\rho^{(N)}$.
The kernel (\ref{ctskernelb}) can be adjusted, by row and column operations, so that for example the 
the diagonal values $\K_{12}(y,y)$ 
are given as the continuous limit of $\K_{12}(y,z)$ as $z \downarrow y$, for example to
 \[
  \left(\begin{array}{ll}
           K^{(c)}_t & -D_2 K^{(c)}_t + \sqrt{\frac{\beta}{\alpha}} K^{(c)}_t \\
           - D_1 K^{(c)}_t + \sqrt{\frac{\beta}{\alpha}} K^{(c)}_t & 
           D_1 D_2 K^{(c)}_t - \sqrt{\frac{\beta}{\alpha}}
            (D_2 K^{(c)}_t + D_1 K^{(c)}_t)
           + \frac{\beta}{\alpha} K^{(c)}_t
           \end{array}\right)
 \]
We hope that this form may be more useful for example when manipulating Fredholm Pfaffians
(as for example in the manipulations for the gap probabilities for the GOE spectrum). 

\textbf{4.2.} With our convention on kernels just as measurable functions, a Poisson rate $\gamma$ process
 can be realised as a Pfaffian point process with kernel $\gamma \J$ where
 $\J(y,z) = 0$ for $y<z$ and $\J_{12}(y,y) = \gamma$. (The same convention would allow
 Poisson processes to be determinantal processes with a purely diagonal kernel. ) 
The kernel in example (b) is connected to Poisson thickening. A locally finite point process 
$X$ can be {\it $\gamma$ thickened} by adding the points of an independent 
Poisson process $Y$ of rate $\gamma$, producing a new point process  $X +Y$. 
If the original point process was Pfaffian with kernel $\K$ then the thickened
process remains Pfaffian with kernel $\K+ \gamma \J$. Indeed, since the points of $X$ and the Poisson process 
never meet, the intensities for the thickened process are given by
\[
\rho^{(N)}_{X+Y} (y_1,\ldots,y_N) = \sum_{J \subseteq \{1,\ldots,N\}} \rho^{|J|}_X(y_j:j \in J) \gamma^{N-|J|}
\]
where $|J|$ is the size of the subset $J$. But the Pfaffian $\Pf(\K+\gamma \J)$ can be expanded by the
Pfaffian sum formula to give exactly this relation.

A locally finite point process $X$ can be {\it $\gamma$ thinned} by removing each point independently with 
probability $\gamma$, producing a new point process which we denote $\Theta_{\gamma}(X)$. 
If $X$ is Pfaffian with kernel $\K$ then the thinned
process  $\Theta_{\gamma}(X)$ remains Pfaffian, with the kernel $\gamma \K$. 

Since the PDE behind both the continuum branching process and the continuum pairwise immigration model
is the same, the heat equation with constant potential, it is not surprising that there is a connection between their fixed time
distribution.  Using thickening and thinning we can state this: let
\begin{eqnarray*}
X_1 &=&  \mbox{the diffusion limit of BCRW with $p=q=\alpha$, $l=r =  2 \epsilon \sqrt{\alpha \beta}$ and $\eta \equiv 1$;} \\
X_2 & = & \mbox{the diffusion limit of ARWPI with $p=q=\alpha$, $m = \epsilon^2 \beta$ and $\eta \equiv 1$;} \\
Y & = & \mbox{a Poisson point process of rate $\frac12 \sqrt{\beta/\alpha}$, independent of $X_2$.}
\end{eqnarray*}
Then, as point processes on $\R$,
\[
\Theta_{1/2} (X_1) \stackrel{\mathcal{D}}{=} X_2 + Y.
\]
The proof is just the verification that the Pfaffian kernels are identical.

There is a nice dynamic coupling argument that connects annihilating Brownian motions with coalescing Brownian motions
(see \cite{TZ:11}), but we do not know a dynamic coupling that explains the above equality of distributions. 

A similar identity also works for carefully chosen Poisson initial conditions, and also for the processes on $\Z$ with suitably chosen
initial conditions (many details are in \cite{thesis}).  

\textbf{(c) Branching and coalescing model with a single initial particle.} 

We take the same parameter choices as in example (a) but choose an initial condition that is a single particle at the origin. The initial conditions for (\ref{epspdeb}) and (\ref{ctskernelb}) change to
\[
K^{(\epsilon)}_0(y,z) = \phi^{z-y} \ind( 0 \not \in [y,z)), \qquad 
K^{(c)}_0(y,z) = e^{\sqrt{\frac{\beta}{\alpha}}(z-y)} \ind(0 \not \in [y,z)).
\]
The explicit solution is 
\begin{equation} \label{oneptkernel}
K^{(c)}_t(y,z) = e^{\sqrt{\frac{\beta}{\alpha}}(z-y)} \left( 1- \psi_t(y) \psi_t(-z)\right) + 
e^{-\sqrt{\frac{\beta}{\alpha}}(z-y)} \psi_t(-y) \psi_t(z), \
\end{equation}
where
\[
\psi_t(x) = 2 \erfc \left( \frac{x-2 \sqrt{\alpha\beta}t}{\sqrt{2 \alpha t}} \right).
\]
As in the discrete setting, the fixed time distribution started from a single site is quite easy to understand. The limit behaviour  
of the leftmost and rightmost particles $\{l_t,r_t\}$, under the parameter scaling we have used, is known to become that of a sticky pair $\{L_t,R_t\}$, that is the solution of
\begin{eqnarray*}
dL_t & = & \ind(L_t \neq R_t) dB^L_t + \ind(L_t = R_t) dB - 2 \sqrt{\alpha \beta} dt, \quad L_0=0, \\
dR_t & = & \ind(L_t \neq R_t) dB^R_t  + \ind(L_t = R_t) dB + 2 \sqrt{\alpha \beta}  dt, \quad R_0=0,
\end{eqnarray*}
where $B^R,B^L,B$ are three independent Brownian motions run at rate $2 \alpha$.
Uniqueness in law holds and $L_t \leq R_t$ for all $t \geq 0$ (see Proposition 2.1 in \cite{SS:08}).
Let $Y$ be an independent Poisson process of rate $2 \sqrt{\frac{\beta}{\alpha}}$. The point process $X^{(c)}_t$
can be constructed as the pair of particles $L_t$ and $R_t$ together with the particles from $Y$ that lie inside $(L_t,R_t)$. Indeed then
\[
\Pp [X_t[y,z) = 0] = \Pp [R_t <y] + \Pp [L_t \geq z] + e^{-2 \sqrt{\frac{\beta}{\alpha}}(z-y)} \Pp [ L_t<y, R_t \geq z].
\]
Comparing this with the formula (\ref{oneptkernel}) for $K^{(c)}_t(y,z) = e^{\sqrt{\frac{\beta}{\alpha}}(z-y)} 
\Pp [X_t[y,z) = 0] $, and using $\psi_t(x) = \Pp [R_t \geq x] = \Pp [L_t <-x]$,
one finds that
\[
\Pp [ L_t<y, R_t \geq z] = \Pp [L_t<y] \Pp [R_t \geq z] - e^{2\sqrt{\frac{\beta}{\alpha}}(z-y)} 
\Pp [ L_t \geq z] \Pp [R_t <y].
\]
This formula, which is straightforward to verify independently, is one way of describing the joint law of $(L_t,R_t)$.

\textbf{(d) Annihilating model with immigration at the origin.} 

We allow immigration only at one site, namely the origin, producing a model we have come to call the Brownian firework.  The immigration rate must be scaled differently to  example (c) in order to see
a non-trivial effect in the continuum limit.  Thus we consider the ARWPI model with parameters $p_x=q_x=\alpha>0$ for all $x$, with $m_0 =  \beta \epsilon^{-1} \geq 0$ and 
$m_x=0$ for all $x \neq 0$, and with zero initial condition. 
From Theorem \ref{theorem:ARWPI} the process $X^{(\epsilon)}_t$ is a Pfaffian point process on $\epsilon \Z$ with kernel 
$\K_t^{(\epsilon)}$ of the form
 \begin{equation} \label{epskernel4}
\K^{(\epsilon)}_t(y,z)
 = \frac{\epsilon}{2} \left(\begin{array}{cc}
           K^{(\epsilon)}_t(y,z) & -D^{(\epsilon)}_2 K^{(\epsilon)}_t(y,z) \\
           - D^{(\epsilon)}_1 K^{(\epsilon)}_t(y,z) & D^{(\epsilon)}_1 D^{(\epsilon)}_2 K^{(\epsilon)}_t(y,z)
           \end{array}\right)
           \quad \mbox{for $y<z$,}    
 \end{equation}
 and $ (\K^{(\epsilon)}_t)_{12}(y,y) = - \frac{\epsilon}{2} D^{(\epsilon)}_2 K^{(\epsilon)}_t(y,y) $, where 
 the function $K^{(\epsilon)}_t(y,z)$ solves, for $y,z \in \epsilon \Z$ with $y \leq z$, 
\[
\left\{ \begin{array}{rcll}
\partial_t K^{(\epsilon)}_t &=&  \alpha \Delta^{(\epsilon)} K^{(\epsilon)}_t - 
2\beta \epsilon^{-1} (\ind(y=0)+\ind(z=0)) K^{(\epsilon)}_t,\\
K^{(\epsilon)}_t(y,y) & = & 1, \\
K^{(\epsilon)}_0(y,z) & = & 1.
\end{array} \right.
\]
The limiting kernel $K^{(c)}_t(y,z)$ solves 
\begin{equation} 
\label{pde4}
\left\{ \begin{array}{rcl}
\partial_t K^{(c)}_t(y,z) &=&  \alpha \Delta  K^{(c)}_t(y,z)  - 2 \beta (\delta_{y=0}+ \delta_{z=0}) K^{(c)}_t(y,z)\\
K^{(c)}_t(y,y) & = & 1 \\
K^{(c)}_0(y,z) & = & 1.
\end{array} \right.
%\nonumber
\end{equation}
This limiting PDE has a distributional potential consisting of delta functions on the $y=0$ and $z=0$ axes. However it has unique bounded continuous mild solutions, which are smooth away from the axes. 
%The right hand derivative $D_2^+ K^{(c)}_t(y,y)$ becomes infinite at $y=0$, and we cannot fit into the bounded hypotheses of %the kernel convergence lemma. 
We first show convergence of $K^{(\epsilon)}_t$. The probabilistic representation of the limiting continuous PDE is 
\[
K^{(c)}_t(y,z) = \E \left[ e^{- \frac{\beta}{\alpha} L^Y_{t \wedge \tau} -  \frac{\beta}{\alpha}L^Z_{t \wedge \tau}} \right]
\]
where $L^Y$ and $L^Z$ are the (semimartingale) local times at zero of two independent Brownian motions $Y$ and $Z$ run at rate $2 \alpha$.
The corresponding formula
\[
K^{(\epsilon)}_t(y,z) = \E \left[ e^{-\frac{\beta}{\alpha}  L^{Y^{(\epsilon)}}_{t \wedge \tau} - \frac{\beta}{\alpha}  L^{(\epsilon)}_{t \wedge \tau}} \right]
\]
holds for random walks $Y^{(\epsilon)},Z^{(\epsilon)}$ on $\epsilon \Z$, jumping right and left each with rate
$\epsilon^{-2} \alpha$, and their local times, for example
\[
L^{Y^{(\epsilon)}}_t = 2 \alpha \epsilon^{-1} \int^t_0 \ind(Y^{(\epsilon)}_s =0) ds.
\]
Weak (and strong) invariance principles for random walks and their local times are a widely studied topic (see survey of results in [3]). 
In our simple concrete setting, the weak convergence  of $(Y^{(\epsilon)},L^{Y^{(\epsilon)}})$ on $D[0,T] \times D[0,T]$
is straightforward to check (characterising the local time via a Tanaka formula). 
Then the weak convergence of $(Y^{(\epsilon)},Z^{(\epsilon)},L^{Y^{(\epsilon)}},L^{Z^{(\epsilon)}})
\to (Y,Z,L^Y,L^Z)$ can be used to check that 
\[
K^{(\epsilon)}_t(y_{\epsilon},z_{\epsilon}) \to K^{(c)}_t(y,z) \quad \mbox{whenever $y_{\epsilon} \to y, z_{\epsilon} \to z$.}
\]
The (complicated) explicit formula below for $K^{(c)}_t(y,z)$ reveals that the intensity
\[
\rho^{(1)}_t(y) = - D_2 K^{(c)}_t(y,y+) \uparrow \infty \quad \mbox{as $y \to 0$.}
\]  
Thus the 
boundedness conditions, on $K$ and its derivatives, in the kernel convergence lemma (Lemma 9 from \cite{GPTZ:15}) can never hold on $\R$. 
However we may consider the kernel on the set $\R \setminus (-\delta,\delta)$ for any $\delta>0$ 
and obtain a limiting point process and limiting kernel on $\R \setminus (-\delta,\delta)$. 
Finally we take $\delta \downarrow 0$ to construct a limiting process and kernel on $\R \setminus \{0\}$
(we use this trick also when doing the further limits $t \to \infty$ and $\beta \to \infty$ below).

%We may however break the equation (\ref{pde4}) into three heat equations, on $\{y<z<0\}$, on $\{y<0,z>0\}$ and on $\{0<y<z\}$,
%each with Dirichlet boundary conditions given by the values of $K^{(c)}$. These equations have $C^2$ solutions 
%on their domains at any fixed time, and the lattice approximations converge suitably (that is they and their derivatives are 
% bounded and converge uniformly)  provided one stays away from the boundaries. 
%We choose $\delta>0$ and consider the restriction of the process on $\R \setminus (-\delta,\delta)$. Then we can consider
% the function $K^{(\epsilon)}_t(y,z)$ as a lattice approximation to the simple heat equation on
% $ (\R \setminus (-\delta,\delta))^2$ and the kernel convergence lemma allows us to construct a limiting Pfaffian point process
% on $\R \setminus (-\delta,\delta)$. Then by consistency we can take $\delta \downarrow 0$ and construct the limiting 
% point process on $\R \setminus \{0\}$. 

The solution to the pde (\ref{pde4}) can be found as follows. The substitution $K^{(c)}_t(y,z) = 1 + \tilde{K}^{(c)}_t(y,z)$ 
yields the equation, on $V_2 = \{(y,z): y<z\}$, 
 \begin{equation} \label{pde7}
\partial_t \tilde{K}^{(c)}_t(y,z) =  \alpha \Delta  \tilde{K}^{(c)}_t(y,z)  - 2 \beta (\delta_{y=0}+ \delta_{z=0}) \tilde{K}^{(c)}_t(y,z)
-  2 \beta (\delta_{y=0}+ \delta_{z=0}), 
\end{equation}
with zero boundary $ \tilde{K}^{(c)}_t(y,y) =0$ and zero initial condition 
$\tilde{K}^{(c)}_0(y,z) =0$.
The Green's kernel for this equation on $V_2$ with zero Dirichlet boundary conditions
can be built out of the one-dimensional kernel $g_t(x,x')$ for the operator $\alpha \Delta - 2 \beta \delta_0$ on $\R$, 
which in turn is given in terms of a standard Brownian motion $B$ and its local time at zero $l^0$ by
\begin{eqnarray*}
g_t(x,x')&=&  \E_x \left[ \delta_{x'}(B_{t/2\alpha}) e^{- (\beta/\alpha) l^0_{t/2 \alpha}} \right]  \\
& = & \sqrt{\frac{\alpha}{\pi t}} e^{- \frac{\alpha (x-x')^2}{t}} - \frac{\beta}{2 \alpha} e^{\frac{(|x|+|x'|)\beta}{\alpha}} e^{\frac{\beta^2 t}{4 \alpha^3}} 
\erfc \left( \frac{\beta t^{1/2}}{2\alpha^{3/2}} + \frac{\alpha^{1/2}(|x|+|x'|)}{t^{1/2}} \right)
\end{eqnarray*}
(for the final formula, where $\erfc(z) = (2/\sqrt{\pi}) \int^{\infty}_{z} \exp(-w^2) dw$, see the handbook \cite{borodin} Equation 1.3.7). 
The Green's kernel for (\ref{pde7}) is given by
$G_t((y,z),(z'z')) = g_t(y,y') g_t(z,z') - g_t(y,z') g_t(z,y')$ and the final solution for the 
inhomogeneous equation (\ref{pde7}) is given using d'Alambert's formula
by
\[
\tilde{K}^{(c)}_t(y,z) =  - 2 \beta \int^t_0 \int_{V_2} G_s((y,z),(y',z')) (\delta_{y'=0} + \delta_{z'=0}) dy' dz' ds.
\] 
This yields a complicated formula for the kernel $K^{(c)}_t(y,z) = 1+ \tilde{K}^{(c)}_t(y,z)$. Two limits lead to more 
attractive explicit kernels. 
The limit as $t \to \infty$ for $K^{(c)}_{t}$, %and its first two derivatives as $t \to \infty$ exist, 
yields the kernel for the steady state of the Brownian firework. Indeed $t \to  K^{(c)}_{\infty}(y,z) $
is decreasing and the limit solves the 
elliptic equation  (still with boundary condition $K^{(c)}_{\infty}(y,y) =1$)
\[
 \alpha \Delta  K^{(c)}_{\infty}(y,z)  =  2 \beta (\delta_{y=0}+ \delta_{z=0}) K^{(c)}_{\infty}(y,z).
\]
The solution is
\[
K^{(c)}_{\infty}(y,z) =1 + \frac{2 \beta}{\pi \alpha} \int^{\infty}_0 e^{-\frac{\beta}{\alpha} s} \left( \arctan \frac{y}{s+|z|} - 
\arctan \frac{z}{s+|y|} \right) ds.
\]
%Again this implies that $X^{(c)}_t \to X_{\infty}^{(c)}$ when restricted to $\R \setminus (-\delta,\delta)$. 
A further limit can be taken as the immigration rate at the origin $\beta$ increases to an infinite rate, 
and this yields the kernel (\ref{firework}) stated in the introduction for the infinite rate Brownian firework.

\end{document}